\documentclass[twocolumn]{autart}

\usepackage{graphicx,subfigure}
\usepackage{color,latexsym,amsfonts,amssymb}
\usepackage{bm}     
\usepackage{amsmath}
\usepackage{booktabs} 
\usepackage{boxedminipage}  
\usepackage{wrapfig} 

\DeclareMathOperator*{\argmin}{argmin}

\newtheorem{theorem}{Theorem}

\begin{document}
\begin{frontmatter}

\title{Mean-Variance Optimization of Discrete Time Discounted Markov Decision Processes \thanksref{footnoteinfo}}
\thanks[footnoteinfo]{This work was supported in part by the National Key
Research and Development Program of China (2016YFB0901900), the
National Natural Science Foundation of China (61573206, 61203039,
U1301254), and the Suzhou-Tsinghua Innovation Leading Action
Project.}
\author[xia]{Li Xia}\ead{xial@tsinghua.edu.cn}
\thanks[xia]{Tel.: +86 10 62793029; fax: +86 10 62796115.}
\address{CFINS, Department of Automation,
TNList, Tsinghua University, Beijing 100084, China}
\begin{abstract}
In this paper, we study a mean-variance optimization problem in an
infinite horizon discrete time discounted Markov decision process
(MDP). The objective is to minimize the variance of system rewards
with the constraint of mean performance. Different from most of
works in the literature which require the mean performance already
achieve optimum, we can let the mean discounted performance equal
any constant. The difficulty of this problem is caused by the
quadratic form of the variance function which makes the variance
minimization problem not a standard MDP. By proving the decomposable
structure of the feasible policy space, we transform this
constrained variance minimization problem to an equivalent
unconstrained MDP under a new discounted criterion and a new reward
function. The difference of the variances of Markov chains under any
two feasible policies is quantified by a difference formula. Based
on the variance difference formula, a policy iteration algorithm is
developed to find the optimal policy. We also prove the optimality
of deterministic policy over the randomized policy generated in the
mean-constrained policy space. Numerical experiments demonstrate the
effectiveness of our approach.
\end{abstract}


\end{frontmatter}

\section{Introduction}\label{section_intro}

The mean-variance optimization is an important problem in stochastic
optimization and its origin can go back to the pioneering work by H.
Markowitz, 1990 Nobel Laureate in Economics, on the modern portfolio
management \cite{Markowitz52}. In financial engineering, the mean
indicates the return of assets and the variance indicates the risk
of assets. The objective of the mean-variance optimization is to
find an optimal policy such that the mean and the variance of system
rewards are optimized in tradeoff and the efficient frontier (a
curve comprised of \emph{Pareto optima}) is obtained.

The mean-variance optimization is first proposed in a static
optimization form in Markowtiz's original paper \cite{Markowitz52},
in which the decision variables are the investment percentage of
securities and the securities returns are described as random
variables with known means and variances. Then, the mean-variance
optimization is further studied in a dynamic optimization form and
Markov decision processes (MDPs) are widely used as an important
analytical model. The difficulty of this problem mainly comes from
the \emph{non-additiveness} of the variance criterion which makes
the principle of \emph{consistent choice} in dynamic programming
invalid \cite{Sobel82,Xia16}. Such invalidness means that the
optimal action selection during $[t+1, \infty)$ may be not optimal
for the action selection during $[t, \infty)$. In the literature,
there are different ways to study the mean-variance optimization.
Many works studied the variance minimization of MDPs after the mean
performance is already maximized \cite{Guo09,Hernandez99,Huang12}.
For such problem, the variance minimization problem can be
transformed to another standard MDP under an equivalent average or
discounted criterion. There are also studies that use the policy
gradient approach to study the mean-variance optimization when the
policy is parameterized \cite{Prashantha13,Tamar12}.

Sobel \cite{Sobel82} gave an early study on the mean-variance
optimization in a discrete time discounted Markov chain, but no
optimization algorithm was presented in that paper. Chung and Sobel
\cite{Chung94,Sobel94} studied the variance minimization problem in
a discrete time Markov chain with the constraint that the long-run
average performance is larger than a given constant. This problem
was transformed to a sequence of linear programming problems, which
may have concerns of computation efficiency since the number of
sequential problems may be large. Guo et al. \cite{Guo15,Huo17}
studied the mean-variance optimization problem in a continuous time
Markov chain with unbounded transition rates and state-action
dependent discount factors, where the performance is accumulated
until a certain state is reached. There are certainly numerous other
excellent works about the mean-variance optimization in the
literature. However, most of the works in the literature either
require a condition of optimal mean performance or reformulate the
problem as variations of mathematical programming. Although linear
programming may be used to study the mean-variance optimization in
some cases, it does not utilize the structure of Markov systems and
the efficiency is not satisfactory. Policy iteration is a classical
approach in dynamic programming and it usually has a high
convergence efficiency. There is little work to study the
mean-variance optimization using policy iteration, at the condition
that the mean performance equals a given value.

In this paper, we study a mean-variance optimization problem in an
infinite horizon discrete time discounted Markov chain. The
objective is to find the optimal policy with the minimal variance of
rewards from the policy set in which the discounted performance
equals a given constant. The motivation of this problem can be
explained with a financial example. People may not always choose an
asset portfolio with the maximal expected return, since a portfolio
with big return usually has big risk (quantified by variance).
People always like to seek a portfolio with minimal risk and
acceptable return. The solution with minimal risk and fixed return
is called \emph{Pareto optimum}. All the Pareto optimal solutions
compose a curve called \emph{Pareto frontier}, or \emph{efficient
frontier} in financial engineering.

The difficulty of such mean-variance optimization problem mainly
comes from two aspects. The first one is the difficulty caused by
the non-additiveness of the variance criterion, which makes the
mean-variance optimization not a standard MDP and policy iteration
is not applicable directly. Another difficulty comes from the fact
that the policy set with a fixed mean performance usually has no
satisfactory structure, such as that described later in
Theorem~\ref{theorem1}. For example, the policy set whose long-run
average performance equals a given constant may not be decomposable
as that in Theorem~\ref{theorem1}. Without such property, dynamic
programming and policy iteration cannot be used for these problems.

In this paper, we use the sensitivity-based optimization theory to
study this nonstandard MDP problem. For the policy set in which the
discounted performance equals a given constant, we prove that this
policy set is decomposable on the action space and the action can be
chosen independently at every state. A difference formula is derived
to quantify the variance difference under any two feasible policies.
The original variance minimization problem with constraints is
transformed to a standard unconstrained MDP under an equivalent
discounted criterion with a new discount factor $\beta^2$ and a new
reward function, where $\beta$ is the discount factor of the
original Markov chain. With this equivalent MDP, we prove the
existence of the optimal policy for this mean-variance optimization
problem. A policy iteration algorithm is developed to find the
optimal policy with the minimal variance. The optimality of
deterministic policy is also proved, compared with randomized
policies generated in the mean-constrained policy space. Finally, we
conduct a numerical experiment to demonstrate the effectiveness of
our approach. The efficient frontier of this numerical example is
also analyzed.

This paper is a continued work compared with our previous papers
\cite{Xia16,Xia17}, which aim to minimize the variance of the
long-run average performance of the Markov chain without considering
the constraint of mean performance. The targeted models in these
papers are different, so are the main results. To the best of our
knowledge, this is the first paper that develops a policy iteration
algorithm to minimize the variance of a discrete time discounted
Markov chain at the condition of any given discounted performance.

\section{Problem Formulation}\label{section_model}

We consider a finite MDP in discrete time. $X_t$ is denoted as the
system state at time $t$, $t=0,1,\cdots$. The state space is finite
and denoted as $\mathcal S=\{1,2,\cdots,S\}$, where $S$ is the size
of the state space. We only consider the deterministic and
stationary policy $d$ which is a mapping from the state space to the
action space. If the current state is $i$, the policy $d$ determines
to choose an action $a$ from a finite action space $\mathcal A(i)$
and a system reward $r(i,a)$ is obtained. The system will transit to
a new state $j$ with a transition probability $p(j|i,a)$ at the next
time epoch, where $i,j \in \mathcal S$ and $a \in \mathcal A(i)$.
Obviously, we have $\sum_{j \in \mathcal S} p(j|i,a) = 1$. Since $d$
is a mapping in the state space, we have $a=d(i)$ and $d(i) \in
\mathcal A(i)$ for all $i \in \mathcal S$. We define the policy
space $\mathcal D$ as the family of all deterministic stationary
policies. For each $d \in \mathcal D$, $\bm P(d)$ is denoted as a
transition probability matrix and its $(i,j)$th element is
$p(j|i,d(i))$, and $\bm r(d)$ is denoted as a column vector and its
$i$th element is $r(i,d(i))$. We assume that the Markov chain is
ergodic for any policy in $\mathcal D$. The discount factor of the
MDP is $\beta$, $0<\beta<1$. For initial state $i$, the mean
discounted performance of the MDP under policy $d$ is defined as
below.
\begin{equation}\label{eq_Jmean}
J(d,i) := \mathbb E^{d}_i\left[ \sum_{t=0}^{\infty} \beta^t
r(X_t,d(X_t)) \right], \quad i \in \mathcal S,
\end{equation}
where $E^{d}_i$ is an expectation operator of the Markov chain at
the condition that the initial state is $i$ and the policy is $d$.
$\bm J(d)$ is an $S$-dimensional column vector and its $i$th element
is $J(d,i)$. The variance of the discounted Markov chain is defined
as below.
\begin{equation}\label{eq_sigma}
\sigma^2(d,i) := \mathbb E^{d}_i \left[ \left(\sum_{t=0}^{\infty}
\beta^t r(X_t,d(X_t))\right) - J(d,i) \right]^2 , \ i \in \mathcal
S.
\end{equation}
We observe that $\sigma^2(d,i)$ quantifies the variance of the
limiting random variable $\sum_{t=0}^{\infty} \beta^t
r(X_t,d(X_t))$. $\bm \sigma^2(d)$ is the variance vector of the
discounted Markov chain and its $i$th element is $\sigma^2(d,i)$.

Denote $\bm \lambda$ as a given mean reward function on $\mathcal
S$. That is, $\bm \lambda$ is an $S$-dimensional column vector and
its $i$th element is denoted as $\lambda(i)$, $i \in \mathcal S$.
The set of all feasible policies with which the mean discounted
performance of the Markov chain equals $\bm \lambda$ is defined as
below.
\begin{equation}
\mathcal D_{\bm \lambda} := \{ \mbox{all } d \in \mathcal D |
J(d,i)=\lambda(i), \mbox{ for all } i \in \mathcal S \}.
\end{equation}
Note that $\mathcal D$ is a deterministic stationary policy set, so
is $\mathcal D_{\bm \lambda}$. In this paper, we do not consider
randomized stationary policies. The optimality of deterministic
policies will be studied in the next section, see
Theorem~\ref{theorem5}. It is easy to see that the policy set
$\mathcal D_{\bm \lambda}$ may be empty if the value of $\bm
\lambda$ is not chosen properly. In this paper, we assume that
$\mathcal D_{\bm \lambda}$ is not empty, which is similar to the
assumption in Markowitz's mean-variance portfolio problem
\cite{Markowitz52,Zhou03}. For a given discounted performance vector
$\bm \lambda$, $\mathcal D_{\bm \lambda}$ may contain more than one
policy. The objective of our mean-variance optimization is to find
an optimal policy from $\mathcal D_{\bm \lambda}$ such that the
variance of the Markov chain is minimized. The mathematical
formulation is written as below.
\begin{equation}\label{eq_problem}
\min\limits_{d \in \mathcal D_{\bm \lambda}} \left\{ \sigma^2(d,i)
\right\}, \quad \mbox{for all } i \in \mathcal S.
\end{equation}
That is, we aim to find an optimal policy among all feasible
policies whose mean discounted performance is equal to a given
constant vector $\bm \lambda$, such that the variance of discounted
rewards is minimized. We denote such a \emph{mean-variance optimal
policy} as $d^*_{\bm \lambda}$. The existence of the solution
$d^*_{\bm \lambda}$ to the problem (\ref{eq_problem}) is not
guaranteed because the minimization in (\ref{eq_problem}) is over
every state $i \in \mathcal S$, i.e., (\ref{eq_problem}) can be
viewed as a multi-objective optimization problem. For a general
multi-objective optimization problem, it is possible that no
solution can dominate all the other solutions on the value of every
dimension of objective function. In the next section, we will
discuss the existence of such optimal policy $d^*_{\bm \lambda}$ and
develop an optimization algorithm to find it. Moreover, we have the
following remarks about this mean-variance optimization problem.

\noindent \textbf{Remark 1.} In the literature on mean-variance
optimization of MDPs, most of the works study the variance
minimization at the constraint that the mean performance is already
maximized. As a comparison in our problem (\ref{eq_problem}), $\bm
\lambda$ can be given as any feasible value.

\noindent \textbf{Remark 2.} If $d^*_{\bm \lambda}$ exists for a
given $\bm \lambda$, then $d^*_{\bm \lambda}$ is also called an
\emph{efficient policy} for the mean-variance optimization problem.
$(\bm \sigma^2(d^*_{\bm \lambda}), \bm \lambda)$ is called an
\emph{efficient point} and all the efficient points with different
$\bm \lambda$'s compose the \emph{efficient frontier} of that
problem.

\section{Main Results}\label{section_mainresult}
We use the sensitivity-based optimization theory to study this
mean-variance optimization problem. This theory was proposed by Dr.
X.-R. Cao and its key idea is the difference formula that quantifies
the performance difference of Markov systems under any two policies
\cite{Cao07}.

Consider a policy $d \in \mathcal D$. For simplicity, we omit $d$ by
default in the following notations $\bm J$, $\bm \sigma^2$, $\bm P$,
$\bm r$. Consider another arbitrary policy $d' \in \mathcal D$, and
the corresponding notations under policy $d'$ are written as $\bm
J'$, ${\bm \sigma^2}'$, $\bm P'$, $\bm r'$. For the mean discounted
performance defined in (\ref{eq_Jmean}), we can rewrite it in a
recursive form as below.
\begin{equation}\label{eq_poisson1}
\bm J = \bm r + \beta \bm P \bm J.
\end{equation}
We further have
\begin{equation}\label{eq_Jinv}
\bm J = (\bm I-\beta\bm P)^{-1} \bm r,
\end{equation}
where $\bm I$ is an $S$-dimensional identity matrix. Note that the
matrix $(\bm I-\beta\bm P)$ is invertible and $(\bm I-\beta\bm
P)^{-1}$ is always nonnegative based on the following observations
\begin{equation}\label{eq_Pinv}
(\bm I-\beta\bm P)^{-1} = \sum_{n=0}^{\infty} \beta^n \bm P^n,
\end{equation}
and $\bm P$ is a nonnegative matrix. We further observe that if the
Markov chain with $\bm P$ is ergodic, all the elements of $(\bm
I-\beta\bm P)^{-1}$ are positive.

For a new policy $d' \in \mathcal D$, we similarly have
\begin{equation}\label{eq_poisson2}
\bm J' = \bm r' + \beta \bm P' \bm J'.
\end{equation}
Subtracting (\ref{eq_poisson1}) from (\ref{eq_poisson2}), we have
\begin{equation}
\bm J' - \bm J = \bm r' - \bm r + \beta(\bm P' - \bm P) \bm J +
\beta \bm P' (\bm J' - \bm J).
\end{equation}
Therefore, we derive the following formula about the difference of
the mean discounted performance of discrete time Markov chains
\begin{equation}\label{eq_diffmean}
\bm J' - \bm J = (\bm I - \beta \bm P')^{-1} \left[ \beta(\bm P' -
\bm P)\bm J + \bm r' - \bm r \right].
\end{equation}
In the above formula, $\bm J$ in the right-hand side is also called
the \emph{performance potential} (or value function) in the
sensitivity-based optimization theory \cite{Cao07}. For a discounted
Markov chain, the performance potential is the same as the mean
discounted performance, but they are different in other general
cases.

\noindent\textbf{Remark 3.} We can derive a policy iteration
algorithm directly based on (\ref{eq_diffmean}). With
(\ref{eq_Pinv}), we see that all the elements of $(\bm I-\beta\bm
P')^{-1}$ are positive for ergodic Markov chains. If we choose $(\bm
P', \bm r')$ to make the vector $\beta \bm P' \bm J + \bm r'$ as
large as possible, then we see $\bm J' - \bm J \geq 0$ and the
policy is improved.

Based on the difference formula (\ref{eq_diffmean}), we study the
structure of the policy set $\mathcal D_{\bm \lambda}$ in which the
mean discounted performance equals a given constant $\bm \lambda$.
The following theorem describes the decomposable structure of
$\mathcal D_{\bm \lambda}$.
\begin{theorem}\label{theorem1}
For any given constant $\bm \lambda$, the policy set $\mathcal
D_{\bm \lambda}$ in which the discounted performance equals $\bm
\lambda$ is the Cartesian product of feasible action space at every
state, i.e.,
\begin{equation}\label{eq_Dlambda}
\mathcal D_{\bm \lambda} = \mathcal A_{\bm \lambda}(1) \times
\mathcal A_{\bm \lambda}(2) \times \cdots \times \mathcal A_{\bm
\lambda}(S),
\end{equation}
where $\mathcal A_{\bm \lambda}(i)$, $i \in \mathcal S$, is defined
as {\scriptsize
\begin{equation}\label{eq_Alambda}
\mathcal A_{\bm \lambda}(i) := \{\mbox{all } a \in \mathcal A(i) |
r(i,a) + \beta \sum_{j \in \mathcal S} p(j|i,a) \lambda(j) =
\lambda(i) \}.
\end{equation}}
\end{theorem}
\vspace{-0.6cm}\noindent \textbf{Proof.} Suppose $d$ is an element
of $\mathcal D_{\bm \lambda}$ and the corresponding transition
probability matrix and reward function are denoted as $\bm P$ and
$\bm r$, respectively. With (\ref{eq_poisson1}), we have
\begin{equation}\label{eq_13}
\bm \lambda = \bm r + \beta \bm P \bm \lambda.
\end{equation}
Consider an arbitrary policy $d' \in \mathcal D$ with $\bm P', \bm
r'$, and the corresponding discounted performance is denoted as $\bm
J'$.

First, we want to prove that if $d' \in \mathcal A_{\bm \lambda}(1)
\times \mathcal A_{\bm \lambda}(2) \times \cdots \times \mathcal
A_{\bm \lambda}(S)$, then $d' \in \mathcal D_{\bm \lambda}$. Since
$d' \in \mathcal A_{\bm \lambda}(1) \times \mathcal A_{\bm
\lambda}(2) \times \cdots \times \mathcal A_{\bm \lambda}(S)$, we
have
\begin{equation}\label{eq_14}
r(i,d'(i)) + \beta \sum_{j \in \mathcal S} p(j|i,d'(i)) \lambda(j) =
\lambda(i), \quad i \in \mathcal S.
\end{equation}
That is, we have
\begin{equation}\label{eq_15}
\bm r' + \beta \bm P' \bm \lambda = \bm \lambda.
\end{equation}
Substituting (\ref{eq_13}) into the above equation, we have
\begin{equation}\label{eq_16}
\bm r' + \beta \bm P' \bm \lambda = \bm r + \beta \bm P \bm \lambda.
\end{equation}
With the difference formula (\ref{eq_diffmean}), we directly have
\begin{equation}\label{eq_17}
\bm J' - \bm \lambda = (\bm I - \beta \bm P')^{-1} \left[ \beta(\bm
P' - \bm P)\bm \lambda + \bm r' - \bm r \right].
\end{equation}
Substituting (\ref{eq_16}) into (\ref{eq_17}), we obtain
\begin{equation}
\bm J' - \bm \lambda = (\bm I - \beta \bm P')^{-1} \bm 0 = \bm 0.
\end{equation}
Therefore, $\bm J' = \bm \lambda$ and $d' \in \mathcal D_{\bm
\lambda}$.

Second, we want to prove that if $d' \in \mathcal D_{\bm \lambda}$,
then $d' \in \mathcal A_{\bm \lambda}(1) \times \mathcal A_{\bm
\lambda}(2) \times \cdots \times \mathcal A_{\bm \lambda}(S)$. Since
$d' \in \mathcal D_{\bm \lambda}$, we know $\bm J' = \bm \lambda$.
Therefore, the difference formula (\ref{eq_17}) becomes
\begin{equation}\label{eq_19}
\beta(\bm P' - \bm P)\bm \lambda + \bm r' - \bm r = (\bm I - \beta
\bm P') (\bm J' - \bm \lambda) = \bm 0.
\end{equation}
Substituting (\ref{eq_13}) into (\ref{eq_19}), we can directly
derive (\ref{eq_15}) and its componentwise form (\ref{eq_14}).
Therefore, we obtain $d' \in \mathcal A_{\bm \lambda}(1) \times
\mathcal A_{\bm \lambda}(2) \times \cdots \times \mathcal A_{\bm
\lambda}(S)$.

Combining the above two results, we have $\mathcal D_{\bm \lambda} =
\mathcal A_{\bm \lambda}(1) \times \mathcal A_{\bm \lambda}(2)
\times \cdots \times \mathcal A_{\bm \lambda}(S)$ and the theorem is
proved. $\hfill \Box$

With Theorem~\ref{theorem1}, we know that $\mathcal D_{\bm \lambda}$
is decomposable as the product of all feasible action spaces
$\mathcal A_{\bm \lambda}(i)$, which indicates that the action
selections at different states in the policy set $\mathcal D_{\bm
\lambda}$ are independent. This property is important for us to
develop a policy iteration algorithm for the problem
(\ref{eq_problem}).

If the performance constraint is quantified under the long-run
average criterion, the decomposable structure of $\mathcal
D_{\lambda}$ may not hold, which makes the policy iteration
inapplicable to such problem. It is partly because the analog of
(\ref{eq_Alambda}) is not valid for the long-run average case. If
the constraint is that the long-run average performance equals the
maximum, the structure of $\mathcal D_{\lambda}$ is still
decomposable and the policy iteration is applicable to the variance
criterion \cite{Guo09,Hernandez99}.

We further study how to minimize the variance of the Markov chain.
Consider a policy $d \in \mathcal D$, the variance is defined as
(\ref{eq_sigma}). The symbol $d$ is also omitted in the following
notations for simplicity. We rewrite (\ref{eq_sigma}) as follows.
\begin{eqnarray}\label{eq_20}
\sigma^2(i) &=& \mathbb E_i \left[ \left(\sum_{t=0}^{\infty} \beta^t
r(X_t)\right) - J(i) \right]^2  \nonumber\\
&=& \mathbb E_i \left[\sum_{t=0}^{\infty} \beta^t r(X_t) \right]^2 -
2 \mathbb E_i \hspace{-0.1cm} \left[\sum_{t=0}^{\infty} \beta^t
r(X_t) \right]\hspace{-0.1cm} J(i) + J^2(i) \nonumber\\
&=& \mathbb E_i \left[ r(X_0) + \sum_{t=1}^{\infty} \beta^t
r(X_t) \right]^2  - J^2(i) \nonumber\\
&=& r^2(i) + 2 \beta r(i) \sum_{j \in \mathcal S} p(j|i) \mathbb E_j
\left[
\sum_{t=0}^{\infty} \beta^t r(X_t) \right] \nonumber\\
&& + \beta^2 \sum_{j \in \mathcal S}p(j|i) \mathbb E_j \left[
\sum_{t=0}^{\infty} \beta^t r(X_t) \right]^2 - J^2(i),
\end{eqnarray}
where the last equality holds using the conditional expectation and
the Markovian property. Recursively substituting (\ref{eq_Jmean})
and (\ref{eq_sigma}) into (\ref{eq_20}), we obtain
\begin{eqnarray}\label{eq_21}
\hspace{-0.1cm} \sigma^2(i) &=& r^2(i) + 2 \beta r(i) \sum_{j \in
\mathcal S} p(j|i) J(j) \nonumber\\
&& + \beta^2 \sum_{j \in \mathcal S}p(j|i) \left[ \sigma^2(j) +
J^2(j) \right] - J^2(i) \nonumber\\
&=& r^2(i) + 2 \beta r(i) \sum_{j \in
\mathcal S} p(j|i) J(j) + \beta^2 \sum_{j \in \mathcal S}p(j|i) J^2(j) \nonumber\\
&& - J^2(i) + \beta^2 \sum_{j \in \mathcal S}p(j|i) \sigma^2(j).
\end{eqnarray}
For the current policy $d$ and state $i$, we define a new reward
function $h(i)$ as below, $i \in \mathcal S$.
{\small\vspace{-0.6cm}
\begin{equation}\label{eq_22}
h(i) = r^2(i) + 2 \beta r(i) \hspace{-0.1cm}\sum_{j \in \mathcal S}
\hspace{-0.05cm} p(j|i) J(j) + \beta^2 \hspace{-0.1cm}\sum_{j \in
\mathcal S}\hspace{-0.05cm}p(j|i) J^2(j) - J^2(i).
\end{equation}}
Substituting the above equation into (\ref{eq_21}), we obtain
\begin{equation}\label{eq_poisson3}
\sigma^2(i) = h(i) + \beta^2 \sum_{j \in \mathcal S}p(j|i)
\sigma^2(j), \quad i \in \mathcal S.
\end{equation}
We observe that the above equation has a recursive form for
$\sigma^2(i)$'s. Comparing (\ref{eq_poisson3}) with
(\ref{eq_poisson1}), we can see that $\sigma^2(i)$ can be viewed the
mean discounted performance of the Markov chain with discount factor
$\beta^2$ and new reward function $h(i)$, $i \in \mathcal S$. With
(\ref{eq_22}), we further define this new reward function in a
vector form as below.
\begin{equation}\label{eq_25a}
\bm h :=  \bm r^2_{\odot} + 2 \beta \bm r \odot (\bm P \bm J) +
\beta^2 \bm P \bm J^2_{\odot} - \bm J^2_{\odot},
\end{equation}
where $\odot$ indicates the \emph{Hadamard product}
(componentwisely), i.e.,
\begin{equation}
\bm J^2_{\odot} := \bm J \odot \bm J :=
(J^2(1),J^2(2),\cdots,J^2(S))^T.
\end{equation}
Based on the above analysis, we directly have the following theorem.
\begin{theorem}\label{theorem2}
The variance of a $\beta$-discounted Markov chain is equivalent to
the mean discounted performance of the same Markov chain with
discount factor $\beta^2$ and reward function $\bm h$ defined in
(\ref{eq_25a}). That is, we have
\begin{equation}\label{eq_26a}
\bm \sigma^2 = \bm h + \beta^2 \bm P \bm \sigma^2.
\end{equation}
\end{theorem}

Note that the new reward function (\ref{eq_22}) can be further
rewritten as below.
\begin{equation}\label{eq_25}
h(i) :=  \sum_{j \in \mathcal S} p(j|i) [r(i)+\beta J(j)]^2 -
J^2(i).
\end{equation}
We can further define a sample path version of the above definition
as below.
\begin{equation}\label{eq_26}
h(X_t) \overset{\mathbb E}{=} [r(X_t)+\beta J(X_{t+1})]^2 -
J^2(X_t),
\end{equation}
where $\overset{\mathbb E}{=}$ means that the equality holds by
expectation. With Theorem~\ref{theorem2} and the above definitions,
we can further rewrite the variance of the Markov chain as below.
\begin{equation}
\sigma^2(i) = \mathbb E_i \left[ \sum_{t=0}^{\infty} \beta^2 h(X_t)
\right], \quad i \in \mathcal S.
\end{equation}

Besides the above sample path version, the value of $\bm \sigma^2$
can also be computed in other ways that are described by the
following theorem.
\begin{theorem}\label{theorem3} \
\begin{itemize}
\item [(a)] With (\ref{eq_25a}), we have
\begin{equation}\label{eq_30}
\bm \sigma^2 = (\bm I - \beta^2 \bm P)^{-1} \bm h.
\end{equation}
\item [(b)] With $\bm f: = \bm r^2_{\odot} + 2 \beta \bm r \odot (\bm P \bm
J)$, we have
\begin{equation}\label{eq_31}
\bm \sigma^2 = (\bm I - \beta^2 \bm P)^{-1} \bm f - \bm J_{\odot}^2.
\end{equation}
\end{itemize}
\end{theorem}
\noindent \textbf{Proof.} (a) Since the matrix $(\bm I-\beta^2 \bm
P)$ is invertible, (\ref{eq_30}) is directly derived from
(\ref{eq_26a}) in Theorem~\ref{theorem2}.

(b) Substituting $\bm h = \bm f + \beta^2 \bm P \bm J^2_{\odot} -
\bm J^2_{\odot}$ into (\ref{eq_30}), we have
\begin{eqnarray}
\bm \sigma^2 &=& (\bm I - \beta^2 \bm P)^{-1}(\bm f + \beta^2 \bm P
\bm J^2_{\odot} - \bm J^2_{\odot}) \nonumber\\
&=& (\bm I - \beta^2 \bm P)^{-1} \bm f + \sum_{n=0}^{\infty}
(\beta^2 \bm P)^n (\beta^2 \bm P \bm J^2_{\odot} - \bm J^2_{\odot})
\nonumber\\
&=& (\bm I - \beta^2 \bm P)^{-1} \bm f - \bm J^2_{\odot}.
\end{eqnarray}
The theorem is proved. $\hfill \Box$

Comparing (\ref{eq_31}) with (\ref{eq_Jinv}), we observe that $\bm
\sigma^2 + \bm J_\odot^2$ equals the mean discounted performance of
the same Markov chain with discount factor $\beta^2$ and reward
function $\bm f$.

Below, we study the variance difference formula of the Markov chain
under any two policies $d, d' \in \mathcal D$. With
Theorem~\ref{theorem2}, we can view the variance as a special form
of the discounted performance. Therefore, we directly apply the
difference formula (\ref{eq_diffmean}) and obtain
\begin{equation}\label{eq_diffvar1}
{\bm \sigma^2}' - \bm \sigma^2 = (\bm I - \beta^2 \bm P')^{-1}
\left[ \beta^2(\bm P' - \bm P)\bm \sigma^2 + \bm h' - \bm h \right],
\end{equation}
where $\bm P'$ and $\bm h'$ are the transition probability matrix
and the equivalent reward function (\ref{eq_25a}) of the Markov
chain under the policy $d'$.

\begin{figure*}[htbp]
\begin{eqnarray}\label{eq_diffvar2}
{\bm \sigma^2}' - \bm \sigma^2 &=& (\bm I - \beta^2 \bm P')^{-1}
\left[ \beta^2(\bm P' - \bm P)\bm \sigma^2 + \bm r'^2_{\odot} + 2
\beta \bm r' \odot (\bm P' \bm \lambda) + \beta^2 \bm P' \bm
\lambda^2_{\odot} - \bm r^2_{\odot} - 2 \beta \bm r \odot (\bm P \bm
\lambda) - \beta^2 \bm P \bm \lambda^2_{\odot}   \right] \nonumber \\
&=& (\bm I - \beta^2 \bm P')^{-1} \left[ \beta^2(\bm P' - \bm P)\bm
(\bm \sigma^2 + \bm \lambda^2_{\odot}) + \bm r'^2_{\odot} + 2 \beta
\bm r'\odot (\bm P' \bm \lambda) - \bm r^2_{\odot} - 2 \beta \bm r
\odot (\bm P \bm \lambda)\right].
\end{eqnarray}
\end{figure*}

\noindent\textbf{Remark 4. } Different from Remark~3, we cannot do
the policy iteration based on (\ref{eq_diffvar1}) because the value
of $\bm h'$ is unknown. We need to pre-compute $\bm J'$ before we
compute the value of $\bm h'$ based on (\ref{eq_25a}). The
computation of $\bm J'$ under every possible policy $d' \in \mathcal
D$ is a brute-force enumeration, which is unacceptable.

Fortunately, our original problem (\ref{eq_problem}) is to find the
optimal policy with the minimal variance from the policy set
$\mathcal D_{\bm \lambda}$. For any $d,d' \in \mathcal D_{\bm
\lambda}$, their mean discounted performances are equal to $\bm
\lambda$. That is, $\bm J = \bm J' = \bm \lambda$, $\forall \ d,d'
\in \mathcal D_{\bm \lambda}$. Applying this fact and (\ref{eq_25a})
to (\ref{eq_diffvar1}), we derive the variance difference formula
(\ref{eq_diffvar2}) under any two policies $\ d,d' \in \mathcal
D_{\bm \lambda}$. Note that (\ref{eq_diffvar2}) is placed on the top
of this page.

To obtain a concise form for (\ref{eq_diffvar2}), we further define
the following column vectors $\bm g$ and $\bm f$ with elements
\begin{equation}\label{eq_g}
g(i) = \sigma^2(i) + \lambda^2(i), \quad i \in \mathcal S,
\end{equation}
\begin{equation}\label{eq_f}
f(i) = r^2(i) + 2 \beta r(i) \sum_{j \in \mathcal S} p(j|i)
\lambda(j), \quad i \in \mathcal S.
\end{equation}
We see that $\bm g$ can be viewed as the performance potential of
the equivalent Markov chain with discount factor $\beta^2$ and cost
function $\bm f$, and (\ref{eq_f}) can be viewed as a special case
of $\bm f$ defined in Theorem~\ref{theorem3}(b) with $\bm J = \bm
\lambda$.

Substituting (\ref{eq_g}) and (\ref{eq_f}) into (\ref{eq_diffvar2}),
we obtain the variance difference formula of Markov chains under any
two policies $d, d' \in \mathcal D_{\bm \lambda}$ as follows.
\begin{equation}\label{eq_diffvar3}
{\bm \sigma^2}' - \bm \sigma^2 = (\bm I - \beta^2 \bm
P')^{-1}\hspace{-0.1cm} \left[ \beta^2(\bm P' - \bm P)\bm g + \bm f'
- \bm f \right].
\end{equation}
Putting the symbol $d$ back to (\ref{eq_g}) and (\ref{eq_f}), we
have {\small \vspace{-0.4cm}
\begin{equation}\label{eq_g2} g(d,i) = \mathbb
E^d_i\left[ \sum_{t=0}^{\infty} \beta^t r(X_t,d(X_t)) \right]^2.
\end{equation}
\begin{equation}\label{eq_f2}
f(i,d(i)) = r^2(i,d(i)) + 2 \beta r(i,d(i)) \sum_{j \in \mathcal S}
p(j|i,d(i)) \lambda(j).
\end{equation}}
From the above equations, we can see that the value of $\bm f'$
(with element $f(i,d'(i))$) under every possible $d' \in \mathcal
D_{\bm \lambda}$ is known since $\bm r'$ (with element $r(i,d'(i))$)
and $\bm P'$ (with element $p(j|i,d'(i))$) are given values. The
difficulty mentioned in Remark~4 is avoided. We can further develop
a policy iteration algorithm to solve this constrained variance
minimization problem (\ref{eq_problem}). The details are described
in Algorithm~1.

\begin{figure*}
\begin{center}
\begin{boxedminipage}{2.1\columnwidth}
\vspace{5pt}
\begin{itemize}

\item For any given feasible constraint value $\bm \lambda$,
arbitrarily choose an initial policy $d^{(0)} \in \mathcal D_{\bm
\lambda}$ and set $k=0$, where $\mathcal D_{\bm \lambda}$ is
determined by $\mathcal D_{\bm \lambda} = \mathcal A_{\bm
\lambda}(1) \times \mathcal A_{\bm \lambda}(2) \times \cdots \times
\mathcal A_{\bm \lambda}(S)$ and $\mathcal A_{\bm \lambda}(i)$ is
determined by (\ref{eq_Alambda}), $i \in \mathcal S$.

\item For the current policy $d^{(k)}$, compute the value of $\bm g(d^{(k)})$
using (\ref{eq_g}) or (\ref{eq_g2}).

\item Generate a new policy $d^{(k+1)}$ using the following policy
improvement
\begin{eqnarray}\label{eq_policyimp}
d^{(k+1)}(i) &:=& \argmin\limits_{a \in \mathcal A_{\bm
\lambda}(i)}\left\{ \beta^2 \sum_{j \in \mathcal S}p(j|i,a)
g(d^{(k)},j) + r^2(i,a) + 2 \beta r(i,a) \sum_{j \in \mathcal S}
p(j|i,a) \lambda(j) \right\}, \quad i \in \mathcal S,
\end{eqnarray}
where we choose $d^{(k+1)}(i) = d^{(k)}(i)$ if possible.
\item If $d^{(k+1)} = d^{(k)}$, stop and output $d^{(k)}$ as the
optimal policy; otherwise, let $k \leftarrow k+1$ and repeat step~2.
\end{itemize}
\vspace{5pt}
\end{boxedminipage}
\vspace{8pt}\\
\textbf{Algorithm 1.} A policy iteration algorithm to solve the
constrained variance minimization problem (\ref{eq_problem}).
\end{center}
\end{figure*}

In the first step of Algorithm~1, we have to compute $\mathcal
A_{\bm \lambda}(i)$ by using (\ref{eq_Alambda}). Since $\bm \lambda$
is given, we can enumerate every action $a \in \mathcal A(i)$ to see
if the equation in (\ref{eq_Alambda}) can hold. The total number of
comparisons used in (\ref{eq_Alambda}) is $\prod_{i \in \mathcal S}
|A(i)|$, which is affordable compared with the value iteration or
policy iteration. The derivation of (\ref{eq_policyimp}) can be
intuitively understood as the policy improvement step in a standard
policy iteration where we aim to equivalently minimize the
discounted performance of the MDP with discount factor $\beta^2$ and
cost function $\bm f$, according to the results of
Theorem~\ref{theorem2} and Theorem~\ref{theorem3}(b).

With the variance difference formula (\ref{eq_diffvar2}) or
(\ref{eq_diffvar3}), we can derive the following theorem about the
existence of the optimal policy and the convergence of Algorithm~1.
\begin{theorem}\label{theorem4}
The optimal policy $d^*_{\bm \lambda}$ for the problem
(\ref{eq_problem}) exists and Algorithm~1 can converge to $d^*_{\bm
\lambda}$.
\end{theorem}
\noindent \textbf{Proof.} First, we prove the convergence of
Algorithm~1. We compare the variance difference of Markov chains
under two policies $d^{(k)}$ and $d^{(k+1)}$ generated in
Algorithm~1. Substituting the policy improvement
(\ref{eq_policyimp}) into (\ref{eq_diffvar2}), we can see that all
the elements of the column vector represented by the square bracket
in the right-hand side of (\ref{eq_diffvar2}) are nonpositive. If
$d^{(k+1)}(i) \neq d^{(k)}(i)$ at some state $i$, then the $i$th
element of that column vector is negative. On the other hand, we
notice the fact that all the elements in the matrix $(\bm I -
\beta^2 \bm P')^{-1}$ are always positive for any ergodic $\bm P'$.
With (\ref{eq_diffvar2}), we have ${\bm \sigma^2}' - {\bm \sigma^2}
< 0$ and the variance is reduced strictly at each iteration. From
Theorem~\ref{theorem1}, we know that $\mathcal D_{\bm \lambda}$ has
a product form of $\mathcal A_{\bm \lambda}(i)$'s and its size is
finite. Therefore, Algorithm~1 can stop within a finite number of
iterations.

Then, we prove the existence of the optimal policy $d^*_{\bm
\lambda}$ and the output of Algorithm~1 is exactly $d^*_{\bm
\lambda}$ when the algorithm stops. We use the contradiction method
to prove it. Assume Algorithm~1 stops at the policy $d$ and $d \neq
d^*_{\bm \lambda}$. Since $d$ is not optimal, from the definition of
the problem (\ref{eq_problem}), we see that there must exist a
policy, say $d' \in \mathcal D_{\bm \lambda}$, such that
${\sigma^2}'(s) < \sigma^2(s)$ at certain state $s \in \mathcal S$.
Therefore, the $s$th element of the vector ${\bm \sigma^2}' - \bm
\sigma^2$ is negative. Furthermore, in (\ref{eq_diffvar3}), since
the element of the matrix $(\bm I - \beta^2 \bm P')^{-1}$ is always
positive, we can derive that some element of the vector represented
by the square bracket in the right-hand side of (\ref{eq_diffvar3})
must be negative. Without loss of generality, we say the $i$th
element of the vector represented by the square bracket in the
right-hand side of (\ref{eq_diffvar3}) is negative. This indicates
that
$$
\beta^2 \bm P'(i,:)\bm g + f'(i) < \beta^2 \bm P(i,:)\bm g + f(i).
$$
Therefore, with (\ref{eq_policyimp}), the above inequality indicates
that the policy $d$ can be further improved by choosing action
$d'(i)$ at state $i$ while remaining the same choices as $d(j)$ for
$j \neq i$. This means that we can still do the policy improvement
(\ref{eq_policyimp}) and Algorithm~1 cannot stop at the current
policy $d$, which contradicts the assumption that Algorithm~1 stops
at $d$. Therefore, the assumption cannot hold. Since we have proved
that Algorithm~1 stops within a finite number of iterations, the
output policy $d$ must be the optimal policy $d^*_{\lambda}$. The
theorem is proved. $\hfill \Box$

Algorithm~1 can be viewed as a special case of the policy iteration
in the traditional MDP theory since we have transformed the original
problem (\ref{eq_problem}) to a standard discounted MDP, as stated
in Theorem~\ref{theorem2}. Algorithm~1 will also have similar
advantages to those of classical policy iteration algorithms, such
as the fast convergence speed. Therefore, Algorithm~1 is an
efficient approach to solve the mean-constrained variance
minimization problem (\ref{eq_problem}).

Based on Theorem~\ref{theorem3}, we can further derive the following
\emph{optimality equation} that the optimal ${\bm \sigma^2}^*$
should satisfy. {\small \vspace{-0.7cm}
\begin{equation}
{\sigma^2}^*(i) = \min\limits_{a \in \mathcal A_{\bm
\lambda}(i)}\left\{ h(i,a) + \beta^2 \sum_{j \in \mathcal S}
p(j|i,a) {\sigma^2}^*(j) \right\},
\end{equation}}
for every state $i \in \mathcal S$. With this optimality equation,
we can also develop a \emph{value iteration} algorithm to solve the
mean-constrained variance minimization problem. The value iteration
algorithm will converge to the optimal value function ${\bm
\sigma^2}^*$ which is exactly the solution to the original problem
(\ref{eq_problem}). The algorithm is analog to that in the classical
MDP theory and we omit the details.

The main results of this paper have been obtained so far. We return
to study a fundamental problem about the optimality of deterministic
policy. All the above results are based on the problem formulation
(\ref{eq_problem}) in Section~\ref{section_model}, where we limit
our optimization in the deterministic policy space $\mathcal D_{\bm
\lambda}$. Below, we extend to study the reward variance of
randomized policies.

\newcounter{TempEqCnt}
\setcounter{TempEqCnt}{\value{equation}} \setcounter{equation}{47}
\begin{figure*}[htbp]\scriptsize
\begin{eqnarray}\label{eq_ith}
\sum_{a \in \mathcal A_{\bm \lambda}(i)} \theta_{i,a} \Bigg\{
\sum_{j \in \mathcal S} [p(j|i,a)-p(j|i,d(i))]\sigma^2(j) + r^2(i,a)
- r^2(i,d(i)) + 2 \beta r(i,a)\sum_{j \in \mathcal
S}p(j|i,a)\lambda(j) - 2 \beta r(i,d(i))\sum_{j \in \mathcal
S}p(j|i,d(i))\lambda(j) \Bigg\}.
\end{eqnarray}
\end{figure*}
\setcounter{equation}{\value{TempEqCnt}}

Here we consider a special category of randomized policies that are
generated from policy space $\mathcal D_{\bm \lambda}$ spanned by
$\mathcal A_{\bm \lambda}(i)$'s. That is, at each state $i$, we
randomly choose actions $a$ from $\mathcal A_{\bm \lambda}(i)$
according to a probability distribution $\theta_{i,a}$, $i \in
\mathcal S$ and $a \in \mathcal A_{\bm \lambda}(i)$. Obviously, we
have $0 \leq \theta_{i,a} \leq 1$ and $\sum_{a \in \mathcal
A_{\lambda}(i)} \theta_{i,a} = 1$, $\forall i \in \mathcal S$. We
denote $\bm \theta$ as a vector composed of elements
$\theta_{i,a}$'s. Different $\bm \theta$ corresponds to different
randomized policy and we denote the randomized policy as $d^{\bm
\theta}$. Based on the sensitivity-based optimization theory, we
obtain Theorem~\ref{theorem5} as follows.

\begin{theorem}\label{theorem5}
For any randomized policy $d^{\bm \theta}$ generated from $\mathcal
D_{\bm \lambda}$, we have $\bm J(d^{\bm \theta}) = \bm \lambda$ and
$\bm \sigma^2(d^{\bm \theta}) \geq \bm \sigma^2(d^*_{\bm \lambda})$.
That is, we do not need to consider $d^{\bm \theta}$'s for problem
(\ref{eq_problem}).
\end{theorem}
\noindent \textbf{Proof.} The transition probability matrix and the
reward function under the randomized policy $d^{\bm \theta}$ are
written as below, respectively.
\begin{equation}\label{eq_Pr}
\begin{array}{l}p^{\bm \theta}(j|i) := \sum_{a \in \mathcal A_{\bm
\lambda}(i)} \theta_{i,a} p(j|i,a).
\\
r^{\bm \theta}(i) := \sum_{a \in \mathcal A_{\bm \lambda}(i)}
\theta_{i,a} r(i,a).
\end{array}
\end{equation}
From definition (\ref{eq_Jmean}), we can derive the following
equation similar to (\ref{eq_Jinv}).
\begin{equation}
\bm J(d^{\bm \theta}) = (\bm I - \beta \bm P^{\bm \theta})^{-1}\bm
r^{\bm \theta},
\end{equation}
where $\bm P^{\bm \theta}$ and $\bm r^{\bm \theta}$ are the
corresponding quantities under the randomized policy $d^{\bm
\theta}$.

First, we compare the mean performance between the randomized policy
$d^{\bm \theta}$ and any deterministic policy $d \in \mathcal D_{\bm
\lambda}$. Similar to the difference formula (\ref{eq_diffmean}), we
can derive
\begin{equation}\label{eq_Jmeanrand}
\bm J(d^{\bm \theta}) - \bm J = (\bm I - \beta \bm P^{\bm
\theta})^{-1} [\beta(\bm P^{\bm \theta} - \bm P) \bm J + \bm r^{\bm
\theta} - \bm r ].
\end{equation}
In the above equation, $(\bm I - \beta \bm P^{\bm \theta})^{-1}$
is a positive matrix and we discuss the element in the square
bracket. The square bracket in (\ref{eq_Jmeanrand}) is an
$S$-dimensional column vector and its $i$th element can be written
as below.
\begin{equation}
\beta\sum_{j \in \mathcal S}(p^{\bm \theta}(j|i) - p(j|i))J(j) +
r^{\bm \theta}(i) - r(i). \nonumber
\end{equation}
Substituting (\ref{eq_Pr}) into the above equation, we obtain
{\small\vspace{-0.8cm}
\begin{equation}
\sum_{a \in \mathcal
A_{\lambda}(i)}\hspace{-0.2cm}\theta_{i,a}\Bigg[\sum_{j \in \mathcal
S}\beta(p(j|i,a) - p(j|i,d(i)))J(j) + r(i,a) - r(i,d(i)) \Bigg].
\nonumber
\end{equation}}
Since $d(i) \in \mathcal A_{\lambda}(i)$ and $J(i) = \lambda(i)$ for
all $i \in \mathcal S$, with (\ref{eq_Alambda}) we can see that the
above equation equals 0. Therefore, we have $\bm J(d^{\bm \theta}) -
\bm J = 0$ and $\bm J(d^{\bm \theta}) = \bm \lambda$.

Then, we compare the reward variance between policy $d^{\bm \theta}$
and $d \in \mathcal D_{\bm \lambda}$. Similar to (\ref{eq_21}), we
can also obtain the variance $\sigma^2(d^{\bm \theta},i)$ under
policy $d^{\bm \theta}$ as below.
\begin{eqnarray}
\sigma^2(d^{\bm \theta},i) = \hspace{-0.3cm}\sum_{a \in \mathcal
A_{\lambda}(i)}\hspace{-0.3cm} \theta_{i,a} \Bigg\{ r^2(i,a) + 2
\beta r(i,a) \sum_{j \in \mathcal
S}p(j|i,a)J(d^{\bm \theta},j) \nonumber\\
+ \beta^2 \sum_{j \in \mathcal S}p(j|i,a)[\sigma^2(d^{\bm \theta},j)
+ J^2(d^{\bm \theta},j)] - J^2(d^{\bm \theta},i)\Bigg\}. \nonumber
\end{eqnarray}
We denote $\bm h^{\bm \theta}$ as a new reward function and its
element $h^{\bm \theta}(i)$, $i \in \mathcal S$, is defined as
below.
\begin{eqnarray}\label{eq_hm}
h^{\bm \theta}(i) &:=& \hspace{-0.3cm} \sum_{a \in \mathcal
A_{\lambda}(i)} \hspace{-0.2cm} \theta_{i,a}\Bigg\{ r^2(i,a) + 2
\beta r(i,a) \sum_{j
\in \mathcal S}p(j|i,a)J(d^{\bm \theta},j)\Bigg\}  \nonumber\\
&+& \beta^2 \sum_{j \in \mathcal S}p^{\bm \theta}(j|i) J^2(d^{\bm
\theta},j) - J^2(d^{\bm \theta},i).
\end{eqnarray}
Therefore, the result in Theorem~\ref{theorem2} also holds for this
randomized policy $d^{\bm \theta}$ and we have
\begin{equation}
\bm \sigma^2(d^{\bm \theta}) = (\bm I - \beta^2 \bm P^{\bm
\theta})^{-1} \bm h^{\bm \theta}.
\end{equation}
Similarly, we can derive the variance difference formula between
policy $d^{\bm \theta}$ and $d$ as below.
\begin{equation}\label{eq_sigmadif}
\bm \sigma^2(d^{\bm \theta}) - {\bm \sigma^2} = (\bm I - \beta^2 \bm
P^{\bm \theta})^{-1} [ \beta^2 (\bm P^{\bm \theta} - \bm P) {\bm
\sigma^2} + \bm h^{\bm \theta} - \bm h].
\end{equation}
Since $(\bm I - \beta^2 \bm P^{\bm \theta})^{-1}$ is a positive
matrix, we study the value of the element of the square bracket in
the above equation. Substituting (\ref{eq_Pr}) and (\ref{eq_hm})
into (\ref{eq_sigmadif}), we can derive (\ref{eq_ith}) to represent
the $i$th element of the square bracket in (\ref{eq_sigmadif}),
where we use the fact $\bm J(d^{\bm \theta}) = \bm J = \bm \lambda$.
Since the value of the large bracket in (\ref{eq_ith}) has no
relation to $d^{\bm \theta}$, we can view it as a given value. With
(\ref{eq_sigmadif}), it is easy to verify that the optimal $\bm
\theta^*$ with the minimal variance must satisfy a necessary
condition: $\theta^*_{i,a} \in \{ 0, 1 \}$, $\forall a \in \mathcal
A_{\bm \lambda}(i), i \in \mathcal S$. That is, the optimal policy
is a deterministic one and we have $\bm \sigma^2(d^{\bm \theta})
\geq \bm \sigma^2(d^*_{\bm \lambda})$. The theorem is proved.
$\hfill \Box$

\renewcommand{\arraystretch}{1.2}
\begin{table*}[htbp]\scriptsize
\caption{The mean and variance of the discounted Markov chain under
every possible policy.}\label{tab1}
\begin{center}
\begin{tabular}{c|*{12}{@{\hspace{4pt}}c}}
\toprule   & $d_1$ & $d_2$ & $d_3$ & $d_4$  & $d_5$ & $d_6$ & $d_7$ & $d_8$ & $d_9$ & $d_{10}$ & $d_{11}$ & $d_{12}$\\
\midrule $\bm J$ & $\begin{pmatrix} 2.5 \\ 4.5 \end{pmatrix}$ & $\begin{pmatrix} 2.2857 \\ 3.4286 \end{pmatrix}$ & $\begin{pmatrix} 2.5 \\ 4.5 \end{pmatrix}$ & $\begin{pmatrix} 2.5 \\ 4.5 \end{pmatrix}$ & $\begin{pmatrix} 2.5 \\ 4.5 \end{pmatrix}$ & $\begin{pmatrix} 2.125 \\ 3.375 \end{pmatrix}$  & $\begin{pmatrix} 2.5 \\ 4.5 \end{pmatrix}$ & $\begin{pmatrix} 2.5 \\ 4.5 \end{pmatrix}$ & $\begin{pmatrix} 2.6172 \\ 4.5234 \end{pmatrix}$ & $\begin{pmatrix} 2.125 \\ 3.375 \end{pmatrix}$ & $\begin{pmatrix} 2.6312 \\ 4.5562 \end{pmatrix}$ & $\begin{pmatrix} 2.6364 \\ 4.5682 \end{pmatrix}$\\
\hline $\bm \sigma^2$ & $\begin{pmatrix} 0.25 \\ 0.25 \end{pmatrix}$ & $\begin{pmatrix} 0.0834 \\ 0.1052 \end{pmatrix}$ & $\begin{pmatrix} 0.25 \\ 0.25 \end{pmatrix}$ & $\begin{pmatrix} 0.2353 \\ 0.0588 \end{pmatrix}$ & $\begin{pmatrix} 0.3222 \\ 0.2556 \end{pmatrix}$ & $\begin{pmatrix} 0.1302 \\ 0.1302 \end{pmatrix}$ & $\begin{pmatrix} 0.3235 \\ 0.2647 \end{pmatrix}$  & $\begin{pmatrix} 0.2963 \\ 0.0741 \end{pmatrix}$ & $\begin{pmatrix} 0.2271 \\ 0.2271 \end{pmatrix}$ & $\begin{pmatrix} 0.1034 \\ 0.1264 \end{pmatrix}$ & $\begin{pmatrix} 0.2316 \\ 0.2316 \end{pmatrix}$ & $\begin{pmatrix} 0.1964 \\ 0.0491 \end{pmatrix}$\\
\bottomrule
\end{tabular}
\end{center}
\end{table*}

\noindent\textbf{Remark 5.} If we consider the randomized policy
generated from the whole policy space $\mathcal D$, the result in
Theorem~\ref{theorem5} may not hold. Because of the quadratic form
of variance functions in this paper, we cannot convert this
mean-constrained variance minimization problem to a linear program
with constraints, which is widely adopted in the literature on
constrained MDPs \cite{Altman99}. The optimality of deterministic
policy and stationary policy is an unsolved problem that needs
further investigation.

\section{Numerical Example}\label{section_numerical}
Consider a discrete time Markov chain with state space $\mathcal
S=\{1,2\}$. The action space is $A(1)=\{1,2,3\}$ and
$A(2)=\{1,2,3,4\}$. The transition probabilities are $p(2|1,a)=a/4$,
$p(1|1,a)=1-a/4$ for $a \in A(1)$, and $p(1|2,a)=a/4$,
$p(2|2,a)=1-a/4$ for $a \in A(2)$. The rewards in state-action pairs
are $r(1,1)=1$, $r(1,2)=\frac{3}{4}$, $r(1,3)=\frac{19}{32}$;
$r(2,1)=\frac{5}{2}$, $r(2,2)=2$, $r(2,3)=3$, $r(2,4)=\frac{13}{4}$.
The discount factor is $\beta = 0.5$. Obviously, the number of total
policies in $\mathcal D$ is $12$ and all the possible policies are
denoted as follows.
\begin{equation}
\begin{array}{l}
d_1 = (1, 1), \ d_2 = (1, 2), \ d_3 = (1, 3), \ d_4 = (1, 4), \\
d_5 = (2, 1), \ d_6 = (2, 2), \ d_7 = (2, 3), \ d_8 = (2, 4), \\
d_9 = (3, 1), \ d_{10} = (3, 2), \ d_{11} = (3, 3), \ d_{12} = (3, 4).\\
\end{array}\nonumber
\end{equation}
Using (\ref{eq_Jinv}) and (\ref{eq_30}), we can compute the value of
$\bm J$ and $\bm \sigma^2$ of the Markov chain under every possible
policy. The computation results are listed in Table~\ref{tab1}.

From Table~\ref{tab1}, we see that there exist some policies under
which the discounted Markov chain has the same mean, even the same
variance. We let $\bm \lambda = (2.5, 4.5)^{T}$ and verify the main
results derived in Section~\ref{section_mainresult}. From
Table~\ref{tab1}, we see that $\mathcal D_{\bm \lambda} = \{d_1,
d_3, d_4, d_5, d_7, d_8 \}$ in which the mean discounted performance
equals $\bm \lambda$. It is easy to verify that $\mathcal A_{\bm
\lambda}(1) = \{1, 2 \}$ and $\mathcal A_{\bm \lambda}(2) = \{1, 3,
4 \}$. Therefore, we have $\mathcal D_{\bm \lambda} = \mathcal
A_{\bm \lambda}(1) \times \mathcal A_{\bm \lambda}(2)$ and the
result of Theorem~\ref{theorem1} is verified in this example. If we
let $\bm \lambda = (2.125, 3.375)^{T}$, we can also verify that
$\mathcal D_{\bm \lambda} = \mathcal A_{\bm \lambda}(1) \times
\mathcal A_{\bm \lambda}(2)$, where $\mathcal D_{\bm \lambda} =
\{d_6, d_{10} \}$, $\mathcal A_{\bm \lambda}(1) = \{2, 3 \}$, and
$\mathcal A_{\bm \lambda}(2) = \{2 \}$.

Below, we verify the policy iteration algorithm to find the optimal
policy for the problem (\ref{eq_problem}) when $\bm \lambda = (2.5,
4.5)^{T}$. We arbitrarily choose an initial policy from the policy
set $\mathcal D_{\bm \lambda}$, say we choose $d^{(0)} = d_5$. We
compute $\bm g(d^{(0)})$ using (\ref{eq_g}) and obtain $\bm
g(d^{(0)}) = (6.5722, 20.5056)^T$. Then we use the policy
improvement (\ref{eq_policyimp}) to generate a new policy as
follows: for state $i=1$, we have $d^{(1)}(1) = \argmin\limits_{a
\in \{1,2\}}\{6.5139, 6.5722 \} = 1$; for state $i=2$, we have
$d^{(1)}(2) = \argmin\limits_{a \in \{1,3,4\}}\{20.5056, 20.5139,$ $
20.3306 \} = 4$. Therefore, $d^{(1)} = (1,4) = d_4$. For $d^{(1)}$,
we repeat the above process and compute $\bm g(d^{(1)}) = (6.4853,
20.3088)^T$. Then we again use (\ref{eq_policyimp}) to generate the
next policy as follows: for state $i=1$, we have $d^{(2)}(1) =
\argmin\limits_{a \in \{1,2\}}\{6.4853, 6.5368 \} = 1$; for state
$i=2$, we have $d^{(2)}(2) = \argmin\limits_{a \in
\{1,3,4\}}\{20.4632, 20.4853, 20.3088 \} = 4$. Therefore, $d^{(2)} =
(1,4) = d_4$. We have $d^{(2)}=d^{(1)}$ and the stopping criterion
is satisfied. Algorithm~1 stops and outputs $d_4$ as the optimal
policy with the minimal variance among the policy set $\mathcal
D_{\bm \lambda}$. From Table~\ref{tab1}, we find that $d_4$ is truly
the optimal policy with the minimal variance $\bm \sigma^2 =
(0.2353, 0.0588)^T$ among the policy set $\mathcal D_{\bm
\lambda}=\{d_1, d_3, d_4, d_5, d_7, d_8 \}$. Therefore, the
convergence of Algorithm~1 is verified in this example and
Algorithm~1 converges to the optimal policy only through 1
iteration.

From Table~\ref{tab1}, we observe another interesting fact that the
policies $d_1$ and $d_3$ have the same mean and variance. To compare
these two policies, we can further consider other high order
performance metrics, such as the $n$th-order bias optimality
\cite{Cao08}. Moreover, we can also investigate the Pareto
optimality of policies from Table~\ref{tab1}. For the objective of
maximizing the mean and minimizing the variance, we can find: $d_2$
dominates $d_{10}$ and $d_{10}$ dominates $d_6$, i.e., $\bm
J(d_{2})>\bm J(d_{10})>\bm J(d_{6})$ and $\bm \sigma^2(d_{2})<\bm
\sigma^2(d_{10})<\bm \sigma^2(d_{6})$; $d_4$ dominates all the other
policies in $\mathcal D_{\bm \lambda}=\{d_1, d_3, d_4, d_5, d_7, d_8
\}$; $d_{12}$ dominates $d_4$, $d_9$, and $d_{11}$. Therefore, the
mean and the variance of policies $\{d_2, d_{12}\}$ comprise the
efficient frontier of this example.

\section{Conclusion}\label{section_conclusion}
In this paper, we study the variance minimization problem of a
discrete time discounted MDP where the mean discounted performance
is equal to a given constant. By transforming this constrained
variance minimization problem to an unconstrained MDP with discount
factor $\beta^2$ and new reward function $\bm h$, we develop a
policy iteration algorithm to efficiently solve such category of
optimization problems. The success of this approach depends on the
decomposable structure of the policy set $\mathcal D_{\bm \lambda}$,
as stated in Theorem~\ref{theorem1}. Such property may not hold if
we consider the long-run average criterion instead of the discounted
criterion. Therefore, the variance minimization problem of MDPs with
a constraint on the long-run average performance is a future
research topic. Another fundamental topic is to study the optimality
of stationary policy for this constrained variance minimization
problem, which may be generally different from the standard
constrained MDPs in the literature. Moreover, how to extend our
results to more general cases, such as continuous time continuous
state Markov processes or finite horizon Markov chains, is another
interesting topic deserving further investigation.

\section*{Acknowledgement}
The author would like to thank the three anonymous reviewers and the
AE for their very constructive comments and suggestions that improve
the paper.

\begin{wrapfigure}{l}{1in}
\vspace{-10pt}
\includegraphics[width=1in,height=1.25in,clip,keepaspectratio]{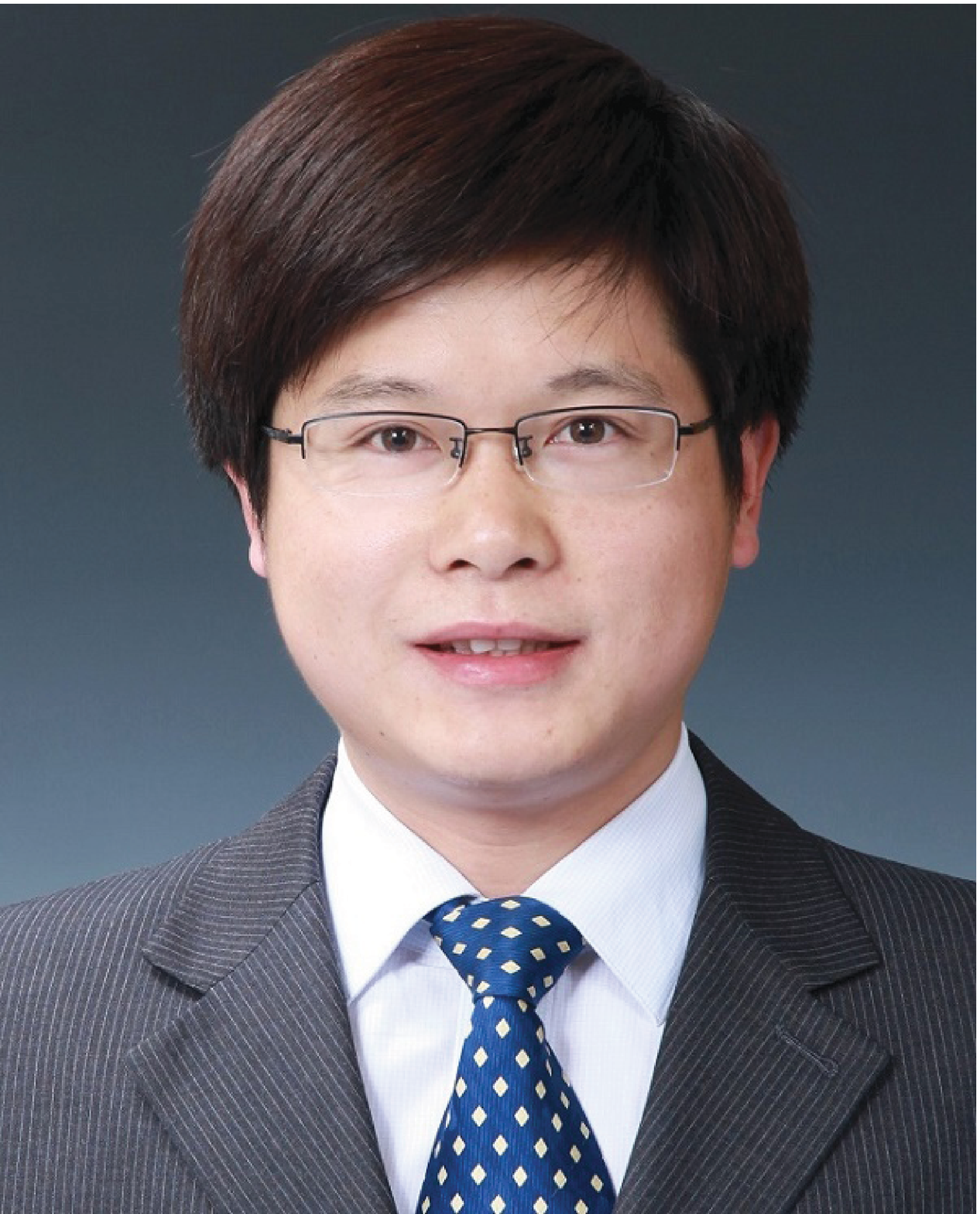}
\vspace{-15pt}
\end{wrapfigure}
\footnotesize \noindent \textbf{Li Xia} is an associate professor in
the Center for Intelligent and Networked Systems (CFINS), Department
of Automation, Tsinghua University, Beijing China. He received the
Bachelor and the Ph.D. degree in Control Theory in 2002 and 2007
respectively, both from Tsinghua University. After graduation, he
worked at IBM Research China as a research staff member (2007-2009)
and at the King Abdullah University of Science and Technology
(KAUST) Saudi Arabia as a postdoctoral research fellow (2009-2011).
Then he returned to Tsinghua University in 2011. He was a visiting
scholar at Stanford University, the Hong Kong University of Science
and Technology, etc. He serves/served as an associate editor and
program committee member of a number of international journals and
conferences. His research interests include the methodology research
in stochastic optimization and learning, queueing theory, Markov
decision processes, reinforcement learning, and the application
research in building energy, energy Internet, industrial Internet,
Internet of things, etc. He is a senior member of IEEE.

\end{document}